%&amstex          
\input amstex\documentstyle{amsppt}  
\pagewidth{12.5cm}\pageheight{19cm}\magnification\magstep1
\topmatter
\title{Lifting involutions in a Weyl group to the torus normalizer}\endtitle
\author G. Lusztig\endauthor
\address{Department of Mathematics, M.I.T., Cambridge, MA 02139}\endaddress
\abstract{Let $N$ be the normalizer of a maximal torus $T$ in a split reductive group over $F_q$
and let $w$ be an involution in the Weyl group $N/T$. We construct explicitly a lifting $n$ of $w$
in $N$ such that the image of $n$ under the Frobenius map is equal to the inverse of $n$.}\endabstract
\thanks{Supported by NSF grant DMS-1566618.}\endthanks
\endtopmatter   
\document

\define\rank{\text{\rm rank}}

\define\dw{\dot w}

\define\ds{\dot s}

\define\mpb{\medpagebreak}

\define\sqc{\sqcup}

\define\qua{\quad}

\define\lb{\linebreak}

\define\op{\oplus}
   
\define\part{\partial}
\define\emp{\emptyset}

\define\ra{\rangle}

\define\m{\mapsto}
\define\do{\dots}
\define\la{\langle}

\define\lra{\leftrightarrow}

\define\sub{\subset}    

\define\T{\times}
\define\ti{\tilde}
\define\nl{\newline}
\redefine\i{^{-1}}
\define\fra{\frac}

\define\ot{\otimes}

\define\Hom{\text{\rm Hom}}

\define\a{\alpha}
\redefine\b{\beta}
\redefine\c{\chi}

\redefine\d{\delta}
\define\e{\epsilon}

\define\io{\iota}
\redefine\o{\omega}

\define\ph{\phi}

\define\s{\sigma}

\define\k{\kappa}

\define\kk{\bold k}

\define\CC{\bold C}

\define\NN{\bold N}

\define\RR{\bold R}

\define\ZZ{\bold Z}
\define\XX{\bold X}
\define\YY{\bold Y}

\define\ce{\Cal E}

\define\cm{\Cal M}
\define\cn{\Cal N}
\define\co{\Cal O}

\define\cz{\Cal Z}

\define\tb{\ti b}

\define\tG{\ti G}

\define\tL{\ti L}

\define\tR{\ti R}
\define\tS{\ti S}
\define\tT{\ti T}
\define\tU{\ti U}

\define\tW{\ti W}
\define\tX{\ti X}
\define\tY{\ti Y}

\define\sha{\sharp}

\define\che{\check}
\define\cha{\che{\a}}

\define\BOU{Bo}
\define\CHEV{Ch}
\define\SQINT{L}
\define\KOS{K}
\define\TIT{T}

\head Introduction\endhead
\subhead 0.1\endsubhead
Let $\kk$ be an algebraically closed field. Let $G$ be a connected reductive
algebraic group over $\kk$.
Let $T$ be a maximal torus of $G$ and let $U$ be the unipotent radical of a Borel subgroup of $G$ containing $T$.
Let $N$ be the normalizer of $T$ in $G$, let $W=N/T$ be the Weyl group, let
$\k:W@>>>N$ be the obvious map. Let $w\m|w|$ be the length function on $W$ and
let $S=\{w\in W;|w|=1\}$. Let $Y=\Hom(\kk^*,T)$. We write the group operation on $Y$ as addition.
 For each $s\in S$ we denote by $\cha_s\in Y$ the corresponding simple coroot; let $L$ be the subgroup of $Y$ generated by $\{\cha_s;s\in S\}$.
Now $W$ acts on $T$ by $w:t\m w(t)=nwn\i$ where $n\in\k\i(w)$; this induces an
action of $W$ on $Y$ and $L$ by $w:y\m y'$ where $y'(z)=w(y(z))$ for $z\in\kk^*$.
We fix a pinning $\{x_s:\kk@>>>G,y_s:\kk@>>>G;s\in S\}$ associated to $T,U$ and we denote by
$w\m\dw$ be the corresponding Tits cross-section \cite{\TIT} of $\k:N@>>>W$.
A {\it halving} of $S$ is a subset $S'$ of $S$ such that $s_1s_2=s_2s_1$ whenever $s_1,s_2$ in $S$ are 
both in $S'$ or both in $S-S'$. Clearly a halving of $S$ exists. 
Let $W_2=\{w\in W;w^2=1\}$. Let $\e=-1\in\kk^*$.
It turns out that, when $w\in W_2$, one can define representatives for $w$
in $\k\i(w)$ other than $\dw$, which in a certain sense are better behaved than $\dw$ (see 0.5).
Namely, for $w\in W_2$, $c\in\kk^*$ and for a halving $S'$ of $S$ we will consider the element
$$n_{w,c,S'}=\dw r_w(c)b_w^{S'}(\e)\in\k\i(w)$$ 
where $r_w\in L$, $b_w\in L/2L$ are given by Theorems 0.2, 0.3 below.
(We then have $r_w(c)\in T$ and $b_w^{S'}(\e)\in T$: if $y\in L$ then $y(\e)\in T$ depends only on the
image of $y$ in $L/2L$ hence $y(\e)\in T$ is defined for any $y\in L/2L$).)

\proclaim{Theorem 0.2} There is a unique map $W_2@>>>L$, $w\m r_w$ such that (i)-(iii) below hold.

(i) $r_1=0$, $r_s=\cha_s$ for any $s\in S$;

(ii) for any $w\in W_2,s\in S$ such that $sw\ne ws$ we have $s(r_w)=r_{sws}$;

(iii) for any $w\in W_2,s\in S$ such that $sw=ws$ we have $r_{sw}=r_w+\cn\cha_s$ where $\cn\in\ZZ$.
\nl
Moreover, in (iii) we have necessarily $\cn\in\{-1,0,1\}$; if in addition $G$ is simply laced we have 
$\cn\in\{-1,1\}$. We have:

(iv) if $w,s$ are as in (iii) and $|sw|>|w|$ then $s(r_w)=r_w$;

(v) if $w\in W_2$ then $w(r_w)=-r_w$.
\endproclaim
A part of the proof of the existence part of the theorem is based on constructing a basis 
consisting of certain positive roots (including the highest root) for the reflection representation of 
$W$ assuming that the longest element is central.
 After I found this basis, I realized  that this basis is the "cascade of roots" that B. Kostant has 
talked about on several occasions. In 2012 he wrote a paper \cite{\KOS} about the cascade. (I thank
D. Vogan for supplying this reference.) The proof of property (iii) is based on a case by case verification.

\proclaim{Theorem 0.3} Let $S'$ be a halving of $S$. There is a unique map
$b=b^{S'}:W_2@>>>L/2L,w\m b_w=b_w^{S'}$ such that (i)-(iii) below hold.

(i) $b_1=0$, $b_s=\cha_s$ for any $s\in S'$, $b_s=0$ for any $s\in S-S'$;

(ii) for any $w\in W_2,s\in S$ such that $sw\ne ws$ we have $s(b_w)=b_{sws}+\cha_s$;

(iii) for any $w\in W_2,s\in S$ such that $sw=ws$ we have $b_{sw}=b_w+l\cha_s$ where $l\in\{0,1\}$;
\nl
Moreover,

(iv) for any $w\in W_2,s\in S$ such that $sw=ws$ we have $s(b_w)=b_w+(\cn+1)\cha_s$ where 
$r_{sw}=r_w+\cn\cha_s$,
$\cn\in\ZZ$;

(v) $b_w(\e)w(b_w(\e))=r_w(\e)\dw^2$, or equivalently $(\dw b_w(\e))^2=r_w(\e)$.
\endproclaim
A part of the proof of this theorem is based on computer calculation.

\subhead 0.4\endsubhead
In this subsection we assume that (i) or (ii) below holds.

(i) $\kk$ is an algebraic closure of a finite field $F_q$ with $q$ elements;

(ii) $\kk=\CC$.

We define $\ph:\kk@>>>\kk$ by $\ph(c)=c^q$ in case (i) and $\ph(c)=\bar c$ (complex conjugation) in case (ii).
In case (i) we assume that $G$ has a 
fixed  $F_q$-rational structure with Frobenius map $\ph:H@>>>H$ such that $\ph(t)=t^q$ for all $t\in T$.

In case (ii) we assume that $G$ has a fixed $\RR$-rational structure so that $G(\RR)$ is the fixed point 
set of an antiholomorphic involution $\ph:G@>>>G$ such that $\ph(y(c))=y(\ph(c))$ for any $y\in Y,c\in\kk^*$.

In both cases we assume that $\ph$ is 
compatible with the fixed pinning of $G$ attached to $T,U$ so that $\ph(\dw)=\dw$ for any $w\in W$. 
In both cases we define $\ph':G@>>>G$ by $\ph'(g)=\ph(g)\i$. 
In case (i), $\ph'$ is a Frobenius map for an $F_q$-rational structure on $G$ which is not in general
compatible with the group structure. In case (ii), $\ph'$ is an antiholomorphic involution of $G$ not
in general compatible with the group structure. Hence $G^{\ph'}=\{g\in G;\ph(g)g=1\}$ is not in general 
a subgroup of $G$.

In both cases we set
$$N^{\ph'}=\{g\in N;\ph(g)g=1\}=N\cap G^{\ph'}.$$ 
Since $\ph(\dw)=\dw$, we see that for $w\in W$, $\k\i(w)\cap N^{\ph'}=\emp$ if $w\in W-W_2$.

We define $\ph':\kk@>>>\kk$ by $\ph'(c)=-\ph(c)$. In case (i) we have
$\kk^{\ph'}=\{x\in\kk;x^q=-x\}$ and in case (ii) we have $\kk^{\ph'}=\{x\in\CC;\bar x+x=0\}$, the set of 
purely imaginary complex numbers. 

Note that for $w\in W_2$, $\dw$ is not necessarily in $N^{\ph'}$. The following result provides some explicit
elements in $\k\i(w)$ which do belong to $N^{\ph'}$.

\proclaim{Theorem 0.5}We assume that we are in the setup of 0.4. Let $w\in W_2$, $c\in\kk^*$ and let $S'$ be a halving of $S$.
We have $\ph'(n_{w,c,S'})=n_{w,\ph'(c),S'}$. Hence if $c\in\k^{\ph'}$ we have $n_{w,c,S'}\in N^{\ph'}$. 
\endproclaim

\subhead 0.6\endsubhead
If $X\sub X'$ are sets and $\io:X'@>>>X'$ satisfies $\io(X)\sub X$ we write $X^\io=\{x\in X;\io(x)=x\}$.

\subhead 0.7\endsubhead
I thank Gongqin Li for help with programming (see 2.4) in GAP using the CHEVIE package \cite{\CHEV}. I also thank 
Meinolf Geck for advice on how to use GAP.

\head 1. The one parameter group $r_w$ attached to an involution $w$ in $W$\endhead
\subhead 1.1\endsubhead
Let $R'$ be a root system in an $\RR$-vector space $\XX'$ of finite dimension; we assume that $R'$ generates 
$\XX$ and that multiplication by $-1$ (viewed as a linear map $\XX'@>>>\XX'$) is contained in the Weyl group $W'$
of $R'$. We assume that we are given a set of positive roots $R'{}^+$ for $R'$. Let $\YY'=\Hom(\XX',\RR)$. 
Let $\la,\ra:\YY'\T\XX'@>>>\RR$ be the obvious pairing. Let 
$\che R'\sub\YY'$ be the set of coroots; let $\a\lra\cha$ be the usual bijection $R'\lra\che R'$. Let 
$\che R'{}^+=\{\cha;\a\in R'{}^+\}$. 
For $\a\in R'$ let $s_\a:\XX'@>>>\XX'$, $s_\a:\YY'@>>>\YY'$ be the reflections defined by $\a$. 

For $\a,\a'$ in $R'{}^+$ we write $\a\le\a'$
if $\a'-\a\in\sum_{\b\in R'{}^+}\RR_{\ge0}\b$. This is a partial order on $R'{}^+$.

Let $\ce_1$ be the set of maximal elements of $R'{}^+$.
For $i\ge2$, let $\ce_i$ is the set of maximal elements of
$$\{\a\in R'{}^+;\la\cha',\a\ra=0 \text{ for any }\a'\in\ce_1\cup\ce_2\cup\do\cup\ce_{i-1}\}.$$
Note that $\ce_1,\ce_2,\do$ are mutually are disjoint. Let 
$$\che\ce_i=\{\cha;\a\in\ce_i\}, \qua \ce=\cup_{i\ge1}\ce_i,\qua \che\ce=\cup_{i\ge1}\che\ce_i.$$
The definition of $\ce,\che\ce$ given above is due to B. Kostant \cite{KOS} who called them {\it cascades}.

From the definition we see that:

(a) {\it if $\a\in\ce,\a'\in\ce,\a\ne\a'$, then $\la\cha',\a\ra=0$.}
\nl
We note the following property.

(b) {\it $\che\ce$ is basis of $\YY'$.}
\nl
For a proof see \cite{\KOS}. Alternatively, we can assume that our root system is irreducible and we can verify
(b) by listing the elements of $\che\ce$ in each case. 
(We denote the simple roots by $\{\a_i;i\in[1,l]\}$ as in \cite{\BOU}.)

Type $A_1$: $\cha_1$.

Type $B_l,l=2n+1\ge3$: $\cha_1+2\cha_2+\do+2\cha_{2n}+\cha_{2n+1}$, 
$\cha_3+2\cha_4+\do+2\cha_{2n}+\cha_{2n+1},\do$, 
$\cha_{2n-1}+2\cha_{2n}+\cha_{2n+1},\cha_1,\cha_3,\do,\cha_{2n+1}$.

Type $B_l,l=2n\ge3$: $\cha_1+2\cha_2+\do+2\cha_{2n-1}+\cha_{2n}$, $\cha_3+2\cha_4+\do+2\cha_{2n-1}+\cha_{2n},\do$,
$\cha_{2n-1}+\cha_{2n},\cha_1,\cha_3,\do,\cha_{2n-1}$.

Type $C_l,l\ge2$: $\cha_1+\cha_2+\do+\cha_{l-1}+\cha_l$, $\cha_2+\do+\cha_{l-1}+\cha_l,\do,\cha_l$.

Type $D_l$, $l=2n\ge4$:  $\cha_1+2\cha_2+2\cha_3+\do+2\cha_{2n-2}+\cha_{2n-1}+\cha_{2n}$, 
$\cha_3+2\cha_4+2\cha_5+\do+2\cha_{2n-2}+\cha_{2n-1}+\cha_{2n}$, 
$\do,\cha_{2n-3}+2\cha_{2n-2}+\cha_{2n-1}+\cha_{2n}$, $\cha_1,\cha_3,\do,\cha_{2n-3},\cha_{2n-1},\cha_{2n}$.

Type $E_7$: $2\cha_1+2\cha_2+3\cha_3+4\cha_4+3\cha_5+2\cha_6+\cha_7$,
$\cha_2+\cha_3+2\cha_4+2\cha_5+2\cha_6+\cha_7$, $\cha_2+\cha_3+2\cha_4+\cha_5$, $\cha_7,\cha_2,\cha_3,\cha_5$.

Type $E_8$: $2\cha_1+3\cha_2+4\cha_3+6\cha_4+5\cha_5+4\cha_6+3\cha_7+2\cha_8$, $2\cha_1+2\cha_2+3\cha_3+4\cha_4+3\cha_5+2\cha_6+\cha_7$,
$\cha_2+\cha_3+2\cha_4+2\cha_5+2\cha_6+\cha_7$, $\cha_2+\cha_3+2\cha_4+\cha_5$, $\cha_7,\cha_2,\cha_3,\cha_5$.

Type $F_4$: $2\cha_1+3\cha_2+2\cha_3+\cha_4$, $\cha_2+\cha_3+\cha_4$, $\cha_2+\cha_3$, $\cha_2$. 

Type $G_2$: $\cha_1+2\a_2,\cha_1$.

\mpb

We note the following result (see \cite{\KOS}).

(c) {\it The reflections $\{s_\b:\YY'@>>>\YY';\b\in\ce\}$ commute with each other and their product
(in any order) is equal to $-1$.}
\nl
Let $\a,\a'$ in $\ce$ be such that $\a\ne\a'$. Using (a) we see that $s_\a(\a')=\a'$. We have also $s_\a(\a)=-\a$. 
Now the result follows from (b).

\subhead 1.2\endsubhead
Let $R',\che R',W'$ be as in 1.1. Let $L'$ be the subgroup of $\YY'$ generated by $\che R'$. We set
$$r=\sum_{\b\in\che\ce}\b\in L'.$$ 
We list the values of $r$ for various types (assumimg that $R'$ is irreducible):

Type $A_1$: $r=\cha_1$.

Type $B_l,l=2n+1\ge3$: $r=2\cha_1+2\cha_2+4\cha_3+4\cha_4+\do+2n\cha_{2n-1}+2n\cha_{2n}+(n+1)\cha_{2n+1}$.

Type $B_l,l=2n\ge4$: $r=2\cha_1+2\cha_2+4\cha_3+4\cha_4+\do+2(n-1)\cha_{2n-3}+2(n-1)\cha_{2n-2}+2n\cha_{2n-1}
+n\cha_{2n}$.

Type $C_l,l\ge2$: $r=\cha_1+2\cha_2+\do+l\cha_l$.     

Type $D_l$, $l=2n\ge4$:  $r=2\cha_1+2\cha_2+4\cha_3+4\cha_4+\do+ (2n-2)\cha_{2n-3}+(2n-2)\cha_{2n-2}
+n\cha_{2n-1}+n\cha_{2n}$.

Type $E_7$: $r=2\cha_1+5\cha_2+6\cha_3+8\cha_4+7\cha_5+4\cha_6+3\cha_7$.

Type $E_8$: $r=4\cha_1+8\cha_2+10\cha_3+14\cha_4+12\cha_5+8\cha_6+6\cha_7+2\cha_8$.

Type $F_4$: $r=2\cha_1+6\cha_2+4\cha_3+2\cha_4$.

Type $G_2$: $r=2\cha_1+2\cha_2$.
\nl
Note that in each case the sum of coefficients of $r$ is equal to $(\sha(R'{}^+)+\rank(R'))/2$.

If $R'$ is irreducible and simply laced, we have $r/2=\sum_{i\in[1,l]}\d_i\o_i$
where $\o_i\in\YY'$ are the fundamental coweights (that is $\la\o_i,\a_j\ra=\d_{ij}$) and $\d_i=\pm1$ are
such that $\d_i+\d_j=0$ when $i,j$ are joined in the Coxeter graph; moreover we have $\d_i=-1$ if in the extended 
(affine) Coxeter graph $i$ is joined with the vertex outside the unextended Coxeter graph.
Another way to state this is that the coefficient of $\cha_i$ in $r$ is equal to half the sum of the coefficients
of the neighbouring $\cha_j$ (that is with $j$ joined with $i$ in the Coxeter graph) plus or minus $1$. For example
in type $E_8$ we have:
$$\align&4=\fra{10}{2}-1, 10=\fra{4+14}{2}+1, 8=\fra{14}{2}+1, 14=\fra{8+10+12}{2}-1, 12=\fra{8+14}{2}+1, \\&
8=\fra{6+12}{2}-1,6=\fra{2+8}{2}+1, 2=\fra{6}{2}-1.\endalign$$
Note that the sign of $\pm1$ in this formula changes when one moves from one $\cha_i$ to a neighbouring one.

\subhead 1.3\endsubhead
In the remainder of this section we place ourselves in the setup of 0.1.
Let $X=\Hom(T,\kk^*)$. We write the group operation in $X$ as addition.
Let $\XX=\RR\ot X,\YY=\RR\ot Y$. Let $\la,\ra:\YY\T\XX@>>>\RR$ be the obvious
nondegenerate bilinear pairing.
The $W$-action on $Y$ in 0.1 induces a linear $W$-action on $\YY$. We define an action of $W$ on $X$ by 
$w:x\m x'$ where $x'(t)=x(w\i(t))$ for $t\in T$. This induces a linear $W$-action on $\XX$.
Let $R\sub X$ be the set of roots; let $\che R\sub Y$ be the set of coroots.
The canonical bijection $R\lra\che R$ is denoted by $\a\lra\cha$.
For any $\a\in R$ we define $s_\a=s_{\cha}:\XX@>>>\XX$ by $x\m x-\la\cha,x\ra\a$ and
$s_\a=s_{\cha}:\YY@>>>\YY$ by $\c\m\c-\la\c,\a\ra\cha$. Then $s_\a=s_{\cha}$ represents the action of an element
of $W$ on $\XX$ and $\YY$ denoted again by $s_\a$ or $s_{\cha}$.
Let $R^+\sub R$ (resp. $\che R^+\sub\che R$) be the set of positive roots (resp corroots) determined by $U$. Let $R^-=R-R^+$, $\che R^-=\che R-\che R^+$.
For $s\in S$ let $\a_s\in R^+$ be the corresponding simple root; 
for any $t\in T$ we have $s(t)=t\cha_s(\a_s(t\i))$; $\cha_s$ has been also considered in 0.1.
Recall that $L$ is the subgroup of $Y$ generated by $\{\cha_s;s\in S\}$.
For any $c\in\kk^*$ we set $c_s=\cha_s(c)\in T$. We have $\ds^2=\e_s$ for $s\in S$. 
Recall that $W_2=\{w\in W;w^2=1\}$. Note that $\dw^2\in T$ for any $w\in W_2$.

\proclaim{Lemma 1.4}Let $w\in W_2$. Then either (i) or (ii) below holds.

(i) There exists $s\in S$ such that $|sw|<|w|$ and $sw\ne ws$.

(ii) There exists a (necessarily unique) subset $J\sub S$ such that $w$ is the longest element in
the subgroup $W_J$ of $W$ generated by $J$; moreover, $w$ is in the centre of $W_J$. For $s\in J$ we 
have $w(\a_s)=-\a_s$. 
\endproclaim
We can assume that $w\ne 1$ and that (i) does not hold for $w$. Let $s_1,s_2,\do,s_k$ in $S$
be such that $w=s_1s_2\do s_k$, $|w|=k$. We have $k\ge1$ and $|s_1w|<|w|$. Since (i) does not hold
we have $s_1w=ws_1$ hence $w=s_2s_3\do s_ks_1$. Thus $|s_2w|<|w|$. Since (i) does not hold
we have $s_2w=ws_2$ hence $w=s_3\do s_{k-1}s_ks_1$. Continuing in this way we see that
$|s_iw|<|w|$ and $s_iw=ws_i$ for $i=1,\do,k$ . We see that the first sentence in (ii) holds with
$J=\{s\in S;s=s_i\text{ for some } i\in[1,k]\}$.

Now let $s\in J$. We have $s=s_\a$ where $\a=\a_s$ and $wsw=s_{w(\a)}$. Since $wsw=s$ we have
$s_{w(\a)}=s_\a$ hence $w(\a)=\pm\a$. Since $|sw|<|w|$ we must have $w(\a)\in R^-$ hence $w(\a)=-\a$.
We see that the second sentence in (ii) holds. The lemma is proved.

\subhead 1.5\endsubhead
For $w\in W_2$ we set $\YY_w=\{y\in\YY;w(y)=-y\}$, $\XX_w=\{x\in\XX;w(x)=-x\}$, $R_w=R\cap\XX_w$,
$\che R_w=\che R\cap\YY_w$, $R_w^+=R^+\cap R_w$, $\che R_w^+=\che R^+\cap\che R_w$.
Note that $\la,\ra$ restricts to a nondegenerate bilinear pairing
$\YY_w\T\XX_w@>>>\RR$ denoted again by $\la,\ra$ and $\a\lra\cha$ restricts to a bijection $R_w\lra\che R_w$.

\proclaim{Lemma 1.6} Let $w\in W_2,s\in S$. Then

(i) $s(\YY_w)=\YY_{sws}$, $s(\XX_w)=\XX_{sws}$, $s(R_w)=R_{sws}$, $s(\che R_w)=\che R_{sws}$;

(ii) if $sw\ne ws$, then $s(R_w^+)=R_{sws}^+$, $s(\che R_w^+)=\che R_{sws}^+$;

(iii) if $sw=ws$ and $|sw|>|w|$, then $s(R_w^+)=R_w^+$, $s(\che R_w^+)=\che R_w^+$;

(iv) if $sw=ws$ and $|sw|<|w|$, then $R_{sw}=\{\a\in R_w;\la\cha_s,\a\ra=0\}$, 
$\che R_{sw}=\{\cha\in\che R_w;\la\cha,\a_s\ra=0\}$, $R_{sw}^+=R_{sw}\cap R_w^+$,
$\che R_{sw}^+=\che R_{sw}\cap\che R_w^+$.
\endproclaim
(i) is immediate. We prove (ii). Let $\a\in R_w^+$; assume that $s(\a)\in R^-$. This implies
that $\a=\a_s$ so that $\a_s\in R_w$ that is $w(\a_s)=-\a_s$ and $w(\cha_s)=-\cha_s$.
For $x\in\XX$ we have
$$\align&ws(x)-sw(x)=w(x-\la\cha_s,x\ra\a_s)-
(w(x)-\la\cha_s,w(x)\ra\a_s)\\&=\la\cha_s,x\ra\a_s+\la w\i\cha_s,x\ra\a_s\\&=
\la\cha_s,x\ra\a_s-\la \cha_s,x\ra\a_s=0.\endalign$$
Thus $ws(x)=sw(x)$ for any $x\in\XX$ so that $sw=ws$ which contradicts our assumption.
We see that $\a\in R_w^+$ implies $s(\a)\in R^+$ hence $s(\a)\in R^+_{sws}$. Thus
$s(R_w^+)\sub R^+_{sws}$. The same argument shows with $w,sws$ interchanged shows that
$s(R_{sws}^+)\sub R^+_w$. It follows that $s(R_w^+)=R^+_{sws}$. Now (ii) follows.

We prove (iii). Let $\a\in R_w^+$; assume that $s(\a)\in R^-$. This implies
that $\a=\a_s$ so that $\a_s\in R_w$ that is $w(\a_s)=-\a_s$. Since $w(\a_s)\in R^-$
we have $|sw|<|w|$ which contradicts our assumption.
We see that $\a\in R_w^+$ implies $s(\a)\in R^+$ hence $s(\a)\in R^+_{sws}=R^+_w$. Thus
$s(R_w^+)\sub R^+_w$. Since $R^+_w$ is finite it follows that $s(R_w^+)=R^+_w$. Now (iii) follows.

We prove (iv). We choose a $W$-invariant positive definite form $(,):\XX\T\XX@>>>\RR$.
Our assumption implies $w(\a_s)=-\a_s$ that is $\a_s\in\XX_w$. Then $\XX_w=\RR \a_s\op\XX'_w$
where $\XX'_w=\{x\in\XX;(x,\a_s)=0\}=\{x\in\XX;\la\cha_s,x\ra=0\}$ and $s$ acts as identity on $\XX'_w$. Since 
$w$ acts as $-1$ on $\XX_w$, $sw$ must act as $-1$ on $\XX'_w$ hence $\XX_{sw}\sub\XX'_w$. Since 
$\dim(\XX_{sw})=\dim(\XX_w)-1=\dim\XX'_w$, it follows that $\XX_{sw}=\XX'_w$. 
We have $R_{sw}=R\cap\XX_{sw}=R\cap\XX'_w$,
$R_{sw}^+=R^+\cap R_{sw}=R^+\cap(R\cap\XX'_w)=R^+\cap\XX'_w$,
$R_{sw}\cap R_w^+=(R\cap\XX_{sw})\cap(R^+\cap\XX_w)=R^+\cap\XX'_w$ hence
$R_{sw}^+=R_{sw}\cap R_w^+$.
Similarly we have $\che R_{sw}=\{\cha\in\che R_w;\la\cha,\a_s\ra=0\}$, $\che R_{sw}^+=\che R_{sw}\cap\che R_w^+$.
This proves (iv).

\proclaim{Lemma 1.7} Let $w\in W_2$.

(a) $R_w$ generates the vector space $\XX_w$ and $\che R_w$ generates the vector space $\YY_w$.

(b) The system $(\YY_w,\XX_w,\la,\ra,\che R_w,R_w)$ is a root system and
$R_w^+$ (resp. $\che R_w^+$) is a set of positive roots (resp. positive coroots) for it.

(c) The longest element of the Weyl group of the root system in (b) acts on $\YY_w$ and on $\XX_w$
as multiplication by $-1$.
\endproclaim
We argue by induction on $|w|$. If $|w|=0$ we have $w=1$ and the lemma is obvious.
Assume now that $|w|>0$. If we can find $s\in S$ such that $|sw|<|w|$, $sw\ne ws$, then by the
induction hypothesis, the lemma is true when $w$ is replaced by $sws$, since $|sws|=|w|-2$. Using
Lemma 1.6 we deduce that the lemma is true for $w$. Using now Lemma 1.4 we see that we can
assume that $w$ is as in 1.4(ii). Let $J\sub S$ be as in 1.4(ii).
Let $\XX'_w$ be the subspace of $\XX_w$ generated by $\{a_s;s\in J\}$. By 1.4(ii) we have
$\XX'_w\sub \XX_w$. As in the proof of 1.4 we can write $w=s_1s_2\do s_k$ with $s_1,s_2,\do,s_k$ in $J$.
Then for any $x\in \XX$ we have
$$\align&wx=s_1s_2\do s_kx=s_2s_3\do s_kx+c_1\a_{s_1}=s_3\do s_kx+c_2\a_{s_2}+c_1\a_{s_1}=\do\\&=
x+c_k\a_{s_k}+\do+c_2\a_{s_2}+c_1\a_{s_1}\endalign$$
with $c_1,c_2,\do,c_k$ in $\RR$. Thus we have 
$(w-1)\XX\sub\XX'_w$. Since $w^2=1$ we have $(w-1)\XX=\XX_w$. Thus $\XX_w\sub\XX'_w$. This proves
the first sentence in (a); the second sentence in (a) is proved in an entirely similar way.
If $\a\in R_w$ (so that $\cha\in \che R_w$) and if $\a'\in R$ then $s_\a(\a')$ is a linear combination
of $\a$ and $\a'$ hence is in $\XX_w$. Since $s_\a(\a')\in R$ we have $s_\a(\a')\in R\cap\XX_w$ that
is $s_\a(\a')\in R_w$. We see that (b) holds. We prove (c). We write again 
$w=s_1s_2\do s_k$ with $s_1,s_2,\do,s_k$ in $J$. We can wiew this as an equality of
endomorphisms of $\XX$ and we restrict it an equality of endomorphisms of $\XX_w$. Each $s_i$ restricts
to an endomorphism of $\XX_w$ which is in the Weyl group of the root system in (b). It follows that
$w$ acts on $\XX_w$ as an element of the Weyl group of the root system in (b) with simple roots
$\{\a_s;s\in J\}$. By 1.4(ii), we have $w(\a_s)=-\a_s$ for any $s\in J$. Thus some element in the
Weyl group of the root system in (b) maps each simple root to its negative. This proves (c).
The lemma is proved.

\mpb

Let $\Pi_w$ be the set of simple roots of $R_w$ such that $\Pi_w\sub R_w^+$.
Let $\che\Pi_w$ be the set of simple coroots of $\che R_w$ such that $\che\Pi_w\sub\che R_w^+$.

\subhead 1.8\endsubhead
Let $w\in W_2$. Let $\che\ce_w$ be the subset of $\che R_w^+$ defined as $\che\ce$ in 1.1 in terms of the root system
$R_w$ in $\XX_w$ (instead of $R'$ in $\XX'$). The definition is applicable in view of 1.7(c). Recall that
$\che\ce_w$ is a basis of $\YY_w$. We define $r_w\in L$ by
$$r_w=\sum_{\b\in\che\ce_w}\b.$$
Note that $r_w$ is a special case of the elements $r$ defined as in 1.2 in terms of $R_w$ instead of $R'$. We 
show:

(a) {\it The reflections $\{s_\b:\YY@>>>\YY;\b\in\che\ce_w\}$ commute with each other and their product
(in any order) is equal to $w$.}
\nl
From 1.1(c) we see that this holds after restriction to $\YY_w$. Since each $s_\b$ and $w$ induces identity on
$\YY/\YY_w$ they must act as $1$ on the orthogonal complement to $\YY_w$ for a $W$-invariant positive definite
inner product on $\YY$. Hence the statements of (a) must hold on $\YY$.

\mpb

We now describe the set $\che\ce_w$ and the elements $r_w$ in the case where $G$ is almost simple and
$w$ is the longest element in $W$ in the case where $w$ is not central in $W$. (The cases where $w$ is central
in $W$ were already described in 1.1, 1.2.) We again denote the simple roots by $\{\a_i;i\in[1,l]\}$ as in
 \cite{\BOU}.

Type $A_l, l=2n\ge2$:
$\cha_1+\cha_2+\do+\cha_{2n-1}+\cha_{2n},\cha_2+\cha_3+\do+\cha_{2n-1}$, $\cha_3+\do+\cha_{2n-2},\do,$
$\cha_{n}+\cha_{n+1}$;

$r_w=\cha_1+2\cha_2+\do+n\cha_n+n\cha_{n+1}+\do+2\cha_{2n-1}+\cha_{2n}$.

Type $A_l, l=2n+1\ge3$:
$\cha_1+\cha_2+\do+\cha_{2n}+\cha_{2n+1},\cha_2+\cha_3+\do+\cha_{2n}$, $\cha_3+\do+\cha_{2n-1},\do,$
$\cha_{n+1}$;

$r_w=\cha_1+2\cha_2+\do+(n+1)\cha_{n+1}+\do+2\cha_{2n}+\cha_{2n+1}$.

Type $D_l, l=2n+1\ge5$:
$\cha_1+2\cha_2+2\cha_3+\do+2\cha_{2n-1}+\cha_{2n}+\cha_{2n+1}$, 
$\cha_3+2\cha_4+2\cha_5+\do+2\cha_{2n-1}+\cha_{2n}+\cha_{2n+1}$, 
$\do,\cha_{2n-3}+2\cha_{2n-2}+2\cha_{2n-1}+\cha_{2n}+\cha_{2n+1}$, $\cha_{2n-1}+\cha_{2n}+\cha_{2n+1}$, 
$\cha_1,\cha_3,\do,\cha_{2n-3},\cha_{2n-1}$;

$r_w=2\cha_1+2\cha_2+4\cha_3+4\cha_4+\do+ (2n-2)\cha_{2n-3}+(2n-2)\cha_{2n-2}
+2n\cha_{2n-1}+n\cha_{2n}+n\cha_{2n+1}$.

Type $E_6$: $\cha_1+2\cha_2+2\cha_3+3\cha_4+2\cha_5+\cha_6$, $\cha_1+\cha_3+\cha_4+\cha_5+\cha_6$,
$\cha_3+\cha_4+\cha_5$, $\cha_4$;

$r_w=2\cha_1+2\cha_2+4\cha_3+6\cha_4+4\cha_5+2\cha_6$.

\mpb

\mpb

From 1.6 and the definitions we deduce the following result.

\proclaim{Lemma 1.9} Let $w\in W_2,s\in S$.

(a) If $sw\ne ws$ then $s(\ce_w)=\ce_{sws}$ and $s(r_w)=r_{sws}$.

(b) If $sw=ws$ and $|sw|>|w|$, then $s(\ce_w)=\ce_w$ and $s(r_w)=r_w$.
\endproclaim

\subhead 1.10\endsubhead
In this subsection we assume that $G$ is almost simple of type $D_l, l\ge4$. Let $\cz=[1,l]$. We can find
a basis $\{e_i;i\in\cz\}$ of $\YY$ with the following properties:

$W$ consists of all automorphisms $w:\YY@>>>\YY$ such that for any $i\in\cz$ we have $w(e_i)=\d_ie_j$ for some 
$j\in\cz$ and some $\d_i\in\{1,-1\}$ and such that $\prod_i\d_i=1$, 
$$\che R^+=\{e_i-e_j;(i,j)\in\cz\T\cz,i<j\}\sqc\{e_i+e_j;(i,j)\in\cz\T\cz,i<j\}.$$
Let $w\in W_2$. Let $P'_w$ be the set of two element subsets $\{i,j\}$ of $\cz$ such that $w(i)=j,w(j)=i$.
Let $P''_w$ be the set of two element subsets $\{i,j\}$ of $\cz$ such that $w(i)=-j,w(j)=-i$.
Let $\cz^+_w=\{i\in\cz;w(i)=i\}$, $\cz^-_w=\{i\in\cz;w(i)=-i\}$. We have
$$\cz=\sqc_{\{i,j\}\in P'_w}\{i,j\}\sqc\sqc_{\{i,j\}\in P''_w}\{i,j\}\sqc\cz^+_w\sqc\cz^-_w.$$
Note that $\sha(\cz^-_w)$ is even. For $w\in W_2$ we have
$$\align&\che\ce_w=\{e_i-e_j;\{i,j\}\in P'_w,i<j\}\sqc\{e_i+e_j;\{i,j\}\in P''_w,i<j\}\sqc \\&
\{e_{i_1}-e_{i_2},e_{i_1}+e_{i_2},e_{i_3}-e_{i_4},
e_{i_3}+e_{i_4},\do,e_{i_{2u-1}}-e_{i_{2u}},e_{i_{2u-1}}+e_{i_{2u}}\}\endalign$$
where $\cz^-_w$ consists of $i_1<i_2<i_3<\do<i_{2u}$. Now let $s\in S$. There are two possibilities:

(i) There exist $a,b$ in $\cz$ such that $b=a+1$, $s(e_a)=e_b,s(e_b)=e_a$, $s(e_z)=e_z$ for
$z\in\cz-\{a,b\}$; moreover, $\cha_s=e_a-e_b$.

(ii) Taking $a=l-1$, $b=l$ we have $s(e_a)=-e_b$, $s(e_b)=-e_a$, $s(e_z)=e_z$ for $z\in\cz-\{a,b\}$;
moreover, $\cha_s=e_a+e_b$.
\nl
We show:

(a) {\it Assume that $w\in W_2,s\in S$ are such that $sw=ws$, $|sw|>|w|$. We have $r_{sw}=r_w\pm\cha_s$.}
\nl
Assume first that $s$ is as in (i). Let $i_1<i_2<i_3<\do<i_{2u}$ be the numbers in $\cz_w^-$. We have either 
$\{a,b\}\in P''_w$ or $\{a,b\}\sub\cz^+_w$. (If $\{a,b\}\in P'_w$ or $\{a,b\}\sub\cz^-_w$ then $|sw|<|w|$.)
If $\{a,b\}\in\cz^+_w$ then $\che\ce_{sw}=\che\ce_w\sqc\{e_a-e_b\}$. Hence $r_{sw}-r_w=\cha_s$.
If $\{a,b\}\in P''_w$ and $i_h<a<b<i_{h+1}$ for some odd $h\in[1,k-1]$ then
$\che\ce_{sw}$ is obtained from $\che\ce_w$ by removing $e_{i_k}-e_{i_{k+1}}$, $e_{i_k}+e_{i_{k+1}}$,
$e_a+e_b$ and by including instead $e_{i_k}-e_a,e_{i_k}+e_a,e_b-e_{i_{k+1}},e_b+e_{i_{k+1}}$. Hence 
$$\align&r_{sw}-r_w=(e_{i_k}-e_a)+(e_{i_k}+e_a)+(e_b-e_{i_{k+1}})+(e_b+e_{i_{k+1}})-\\&
(e_{i_k}-e_{i_{k+1}})-(e_{i_k}+e_{i_{k+1}})-(e_a+e_b)=e_b-e_a=-\cha_s.\endalign$$
If $\{a,b\}\in P''_w$ and there is no odd $h\in[1,k-1]$ such that $i_h<a<b<i_{h+1}$ then
$\che\ce_{sw}=\che\ce_w\sqc\{e_a-e_b\}$. Hence $r_{sw}-r_w=\cha_s$.
Next we assume that $s$ is as in (ii). We have $\che\ce_{sw}=\che\ce_w\sqc\{e_a+e_b\}$. Hence $r_{sw}-r_w=\cha_s$.
This completes the proof of (a).

\subhead 1.11\endsubhead
In this subsection we assume that $G$ is simple of type $E_8$. Let $w$ be the longest element of $W$.
We denote the simple roots by $\{\a_i;i\in[1,8]\}$ as in \cite{\BOU} and we write $s_i$ instead of $s_{\a_i}$.
For $i\in[1,8]$ we have $s_iw\in W_2$ and $R_{s_iw}=\{\a\in R;\la \cha_i,\a\ra=0\}$. From this $\che\Pi_{s_iw}$ is
easily determined in each case:
$$\che\Pi_{s_1w}=\{\cha_1+2\cha_2+2\cha_4+\cha_3+\cha_5,\cha_2,\cha_4,\cha_5,\cha_6,\cha_7,\cha_8\}.$$
$$\che\Pi_{s_2w}=\{\cha_1, \cha_2+2\cha_4+\cha_3+\cha_5,\cha_3,\cha_5,\cha_6,\cha_7,\cha_8\}.$$
$$\che\Pi_{s_3w}=\{\cha_1+\cha_3+\cha_4,\cha_2,\cha_3+2\cha_4+\cha_2+\cha_5,\cha_5,\cha_6,\cha_7,\cha_8\}.$$
$$\che\Pi_{s_4w}=\{\cha_1,\cha_2+\cha_4+\cha_3, \cha_3+\cha_4+\cha_5,\cha_6,\cha_2+\cha_4+\cha_5,\cha_7,\cha_8\}.$$
$$\che\Pi_{s_6w}=\{\cha_1,\cha_2,\cha_3,\cha_4,  \cha_5+\cha_6+\cha_7,\cha_2+\cha_3+2\cha_4+2\cha_5+\cha_6,\cha_8\}.
$$
$$\che\Pi_{s_7w}=\{\cha_1,\cha_2,\cha_3,\cha_4,\cha_5,\cha_6+\cha_7+\cha_8,\cha_2+\cha_3+2\cha_4+2\cha_5+2\cha_6
+\cha_7\}.$$
$$\che\Pi_{s_8w}=\{\cha_1,\cha_2,\cha_3,\cha_4,\cha_5,\cha_6,\cha_2+\cha_3+2\cha_4+2\cha_5+2\cha_6+2\cha_7
+\cha_8\}.$$
This is the set of simple roots of a root system of type $E_7$ hence $r_{s_iw}$ is given by substituting
the simple roots in the formula for $r$ in type $E_7$ given in 1.2 by the roots in $\che\Pi_{s_1w}$. We find
the same result as for $r_w$ (given by $r$ in type $E_8$ in 1.2) plus a multiple of $\cha_i$.
More precisely:
$$r_{s_1w}=r_w-\cha_1,\qua r_{s_2w}=r_w+\cha_2,\qua r_{s_3w}=r_w+\cha_3,$$
$$r_{s_4w}=r_w-\cha_4,\qua r_{s_5w}=r_w+\cha_5,\qua r_{s_6w}=r_w-\cha_6,$$
$$r_{s_7w}=r_w+\cha_7,\qua r_{s_8w}=r_w-\cha_8.$$

\subhead 1.12\endsubhead
In this subsection we preserve the setup and notation of 1.11.
Let $w'$ be the longest element in the standard parabolic subgroup of type $E_7$ of $W$.
For $i\in[1,7]$ we have $s_iw'\in W_2$ and $R_{s_iw'}=\{\a\in R_{w'};\la \cha_i,\a\ra=0\}$. From this 
$\che\Pi_{s_iw'}$ is easily determined in each case:
$$\che\Pi_{s_1w'}=\{\cha_1+2\cha_3+2\cha_4+\cha_2+\cha_5,\cha_2,\cha_4,\cha_5,\cha_6,\cha_7\}.$$
$$\che\Pi_{s_2w'}=\{\cha_1,\cha_3,\cha_3+2\cha_4+\cha_2+\cha_5,\cha_5,\cha_6,\cha_7\}.$$
$$\che\Pi_{s_3w'}=\{\cha_1+\cha_3+\cha_4,\cha_2,\cha_3+2\cha_4+\cha_2+\cha_5,\cha_5,\cha_6,\cha_7\}.$$
$$\che\Pi_{s_4w'}=\{\cha_1,\cha_3+\cha_4+\cha_5, \cha_4+\cha_2+\cha_3,\cha_4+\cha_2+\cha_5,\cha_6,\cha_7\}.$$
$$\che\Pi_{s_5w'}=\{\cha_1,\cha_2,\cha_3,\cha_4+\cha_5+\cha_6, \cha_5+2\cha_4+\cha_2+\cha_3,\cha_7\}.$$
$$\che\Pi_{s_6w'}=\{\cha_1,\cha_2,\cha_3,\cha_4, \cha_5+\cha_6+\cha_7,\cha_6+2\cha_5+2\cha_4+\cha_2+\cha_3\}.$$
$$\che\Pi_{s_7w'}=\{\cha_1,\cha_2,\cha_3,\cha_4,\cha_5,  \cha_7+2\cha_6+2\cha_5+2\cha_4+\cha_2+\cha_3\}.$$
This is the set of simple roots of a root system of type $D_6$ hence $r_{s_iw'}$ is given by substituting
the simple roots in the formula for $r$ in type $D_6$ given in 1.2 by the roots in $\che\Pi_{s_1w'}$. We find
the same result as for $r_{w'}$ (given by $r$ in type $E_7$ in 1.2) plus a multiple of $\cha_i$.
More precisely:
$$r_{s_1w'}=r_{w'}-\cha_1,\qua r_{s_2w'}=r_{w'}+\cha_2,\qua r_{s_3w'}=r_{w'}+\cha_3,\qua
r_{s_4w'}=r_{w'}-\cha_4,$$
$$r_{s_5w'}=r_{w'}+\cha_5,\qua r_{s_6w'}=r_{w'}-\cha_6,\qua r_{s_7w'}=r_{w'}+\cha_7.$$

\subhead 1.13\endsubhead
In this subsection we preserve the setup and notation of 1.11. We show:

(a) {\it Let $z\in W_2,s\in S$ be such that $sz=zs$. We have $r_{sz}=r_z+\cn\cha_s$ where 
$\cn\in\{-1,1\}$.}
\nl
By interchanging if necessary $z,sz$, we can assume that $|z|>|sz|$.
By 1.4 we can find a sequence $s_1,s_2,\do,s_k$ in $S$ (with $k\ge0$) such that
$|z|>|s_1zs_1|>|s_2s_1zss_1s_2|>\do>|s_k\do s_2s_1zs_1s_2\do s_k|$ and
$z':=s_k\do s_2s_1zs_1s_2\do s_k$ is the longest element of a standard parabolic subgroup $W_J$ of $W$ such
that $z'$ is in the centre of $W_J$. Let $\s=s_k\do s_2s_1\in W$.
Applying 1.6(ii) repeatedly we see that
$\che R_{s_k\do s_2s_1zs_1s_2\do s_k}^+=s_k(\che R_{s_{k-1}\do s_2s_1zs_1s_2\do s_{k-1}}^+)$,
$\do,$ $\che R_{s_2s_1zs_1s_2}^+=s_2(\che R_{s_1zs_1}^+)$, $\che R_{s_1zs_1}^+=s_1(\che R_z^+)$.
It follows that 
$\che R_{s_k\do s_2s_1zs_1s_2\do s_k}^+=s_k\do s_2s_1(\che R_z^+)$ that is
$\che R_{z'}^+=s_k\do s_2s_1(\che R_z^+)=\s(\che R_z^+)$. This implies that

(b) $\che \Pi_{z'}=\s(\che\Pi_z)$. 
\nl
From our assumption we have $z(\cha_s)=-\cha_s$. Thus $\cha_s\in \che R_z^+$. Since $\cha_s$ is a simple coroot
in $\che R$ we necessarily have $\cha_s\in\che\Pi_z$. Using (b) we deduce that $\s(\cha_s)\in\che\Pi_{z'}$.
From the definition of $z'$ we see that $\che\Pi_{z'}$ consists of the simple coroots of $\che R$ such that
the corresponding simple reflections are in $W_J$. Thus we have $\s(\cha_s)=\cha_{s'}$ where $s'\in S\cap W_J$. 
It follows that $\s s\s\i=s'$, $\s(\a_s)=\a_{s'}$. Note that $s'z'=z's'$ (since $z'$ is in the centre of $W_J$)
and $|s'z'|<|z'|$ (since $z'$ is the longest element of $W_J$).
From 1.6(iv) we see that $\che R_{sz}^+=\{\cha\in R_z^+;\la\cha,\a_s\ra=0\}$,
$\che R_{s'z;}^+=\{\cha\in R_{z'}^+;\la\cha,\a_{s'}\ra=0\}$.
If $\la\cha,\a_s\ra=0$, then $\la\s(\cha),\s(\a_s)\ra=0$ hence $\la \s(\cha),\a_{s'}\ra=0$.
Since $\s(R_z^+)=R_{z'}^+$ it follows that
$\s(\{\cha\in R_z^+;\la\cha,\a_s\ra=0\})=\{\cha\in R_{z'}^+;\la\cha,\a_{s'}\ra=0\}$ that is,
$\s(\che R_{sz}^+)=\che R_{s'z'}^+$. This implies 
$\s(\che\Pi_{sz})=\che\Pi_{s'z'}$. Using the definitions we deduce that $\s(\ce_{sz})=\ce_{s'z'}$ hence
$\s(r_{sz})=r_{s'z'}$. Similarly we have $\s(r_z)=r_{z'}$. 
Hence if (a) holds for $z',s'$ that is $r_{s'z'}=r_{z'}+\cn\cha_{s'}$ where $\cn\in\{-1,1\}$ then
$r_{sz}-r_z=\s\i(\cn\cha_{s'})=\cn\cha_s$ so that (a) holds for $z,s$.
Thus it is enough to prove (a) assuming in addition that $z$ is the longest element of a standard parabolic 
subgroup $W_J$ of $W$ such that $z$ is in the centre of $W_J$. 
If $W_J=W$, (a) follows from 1.11. If $W_J$ is of type $E_7$, (a) follows from 1.12.
If $W_J$ is of type other than $E_8,E_7$, then it is of type $A_1\T A_1\T\do$ or of type $D_l\T A_1\T A_1\T\do$  
(with $l\in\{4,6\}$). If $s$ belongs to the $A_1\T A_1\T\do$-factor, the result is trivial. If $s$ belongs to
the $D_l$-factor, the result follows from 1.10. This completes the proof of (a).

\subhead 1.14\endsubhead
In this subsection we assume that $G$ is almost simple, simply connected, simply laced and that we are given an 
automorphism $\io:G@>>>G$ such that $\io(T)=T$, $\io(U)=U$ and that
for any $s\in S,c\in\kk^*$ we have $x_{\io(s)}(c)=\io(x_s(c))$,
$y_{\io(s)}(c)=\io(y_s(c))$. Then $\io$ induces
an automorphism of $W$ and automorphisms of $X$ and $Y$ leaving stable $R$, $\che R$, $S$; these are denoted
again by $\io$. 
We also assume that if $s,s'$ in $S$ are in the same $\io$-orbit then $ss'=s's$.
Let $\tG=G^\io$, a connected simply connected algebraic group. Now $\tT=T^\io$ is a maximal torus of $\tG$ and 
$\tU=U^\io$ is the unipotent radical of a Borel subgroup of $\tG$. Let $\tW$ be the Weyl group of $\tG$ with
respect to $\tT$. We can identify $\tW=W^\io$. Let $w\m|w|_\io$ be the length function on $\tW$.
Let $\tS=\{w\in\tW;|w|_\io=1\}$. Now $\tS$ consists of the elements $\s=\prod_ss$ where $s$ runs over an 
$\io$-orbit in $S$. Let $\tX=\Hom(\tT,\kk^*)$ (a quotient of $X$) and let $\tY=\Hom(\kk^*,\tT)$ (a subgroup of
$Y$); we have $\tY=Y^\io$.
Let $\tR$ (resp, $\che{\tR}$) be the set of roots (resp. coroots) of $\tG$ with respect to $\tT$.
Now $\tR$ consists of the images of roots of $G$ under $X@>>>\tX$ and $\che{\tR}$ consists of
the elements of $Y$ which are sums of coroots in an $\io$-orbit on $\che R$.
If $\s\in\tS$ corresponds to a $\io$-orbit $\co$ in $S$ then the simple root $\a_\s$ of $\tG$ corresponding to 
$\s$ is the restriction to $\tX$ of $\a_s$ for any $s\in\co$; the simple coroot of $\tG$ corresponding to $\s$ is 
$\cha_\s=\sum_{s\in\co}\cha_s\in\tY$. Let $\tW_2=W_2\cap\tW$.
Let $\{\ti r_w;w\in \tW_2\}$ be the elements of $\tY$ defined like $\{r_w;w\in W_2\}$ (see 1.8) in terms of
$\tG$ instead of $G$. We show:

(a) {\it For $w\in\tW_2$ we have $\ti r_w=r_w$.}
\nl
We argue by induction on $|w|_\io$. If $|w|_\io=0$ we have $w=1$ and the result is obvious. Assume now that
$|w|_\io\ge1$. Assume also that we can find $\s\in\tS$ such that $|\s w|_\io<|w|_\io$ and $\s w\ne w\s$.
 We write $\co=\{s_1,\do,s_k\}\sub S$, $\s=s_1\do s_k$.
Then for some $i\in[1,k]$ we have $s_iw\ne ws_i$. Hence for all $i\in[1,k]$ we have $s_iw\ne ws_i$. 
Hence we have $s_1w\ne ws_1, s_2s_1w\ne ws_1s_2,\do,s_k\do s_1w\ne ws_1\do s_k$. By 1.9(a) for $G$ and $\tG$
we have $\ti r_w=\s\ti t_{sws}\s$ and
$r_w=s_1(r_{s_1ws_1})=s_1s_2(r_{s_2s_1ws_1s_2})=\do=s_1\do s_k(r_{s_k\do s_1ws_1\do s_k})$ so that 
$r_w=\s r_{\s w\s}s$. By the induction hypothesis we have $\ti r_{\s w\s}=r_{\s w\s}$ hence $\ti r_w=r_w$.
Next we assume that there is no $\s\in\tS$ such that $|\s w|_\io<|w|_\io$ and $\s w\ne w\s$. Then, by 1.4 for
$\tG$ we can find a standard parabolic subgroup $\tW'$ of $\tW$ such that $w$ is the longest element of $\tW'$
and $w$ is central in $\tW'$. In this case the equality $\ti r_w=r_w$ follows by comparing the formulas
in 1.2 with those in 1.10. This completes the proof of (a).

We return to the general case.

\proclaim{Lemma 1.15} Let $w\in W_2,s\in S$ be such that $sw=ws$. We have 

(a) $r_{sw}=r_w+\cn\cha_s$ where $\cn\in\{-1,0,1\}$. 
\nl
If in addition $G$ is simply laced then $\cn\in\{-1,1\}$.
\endproclaim
If the result holds when $|w|>|sw|$ then it also holds when $|w|<|sw|$ (by interchanging $w,sw$); thus we can
assume that $|w|>|sw|$. We can assume that $G$ is almost simple.
If $G$ is of type $D_l$, $l\ge4$, the result follows from 1.10. If $G$ is of type $A_l$, the result follows
from the corresponding result for a group of type $D_{l'}$ with $l<l'\ge4$. If $G$ is of type $E_8$, the
result follows from 1.13(a). If $G$ is of type $E_7$ or $E_6$ the result follows
from the corresponding result for a group of type $E_8$. Thus (a) holds when $G$ is simply laced.

Let $G,\io,\tG,\tT,\tW,\tS,\tY,||_\io,\tW_2$ be as in the proof of 1.14(a). 
To complete the proof it is enough to show that (a) holds when $G$ is replaced by $\tG$.
Let $\{\ti r_w;w\in \tW_2\}$ be the elements of $\tY$ defined like $\{r_w;w\in W_2\}$ in terms of
$\tG$ instead of $G$. We must show that $\{\ti r_w;w\in \tW_2\}$ satisfy conditions like (a).
By 1.14(a) we have $\ti r_w=r_w$ for $w\in\tW_2$.

Now let $w\in\tW_2$, $\s\in\tS$ be such that $\s w=w\s$.
We write $\co=\{s_1,\do,s_k\}\sub S$, $\s=s_1\do s_k$
where $\io$ permutes $s_1,s_2,\do,s_k$ cyclically: $s_1\m s_2\m\do s_k\m s_1$. (Note that $k\le 3$.)
We have $w(\cha_\s)=\cha_\s$ hence $w(\cha_{s_1}+\do+\cha_{s_k})=\cha_{s_1}+\do+\cha_{s_k}$.
If $w(\cha_{s_i})\in-\che R^+$ for some $i\in[1,k]$ then the same is true for any $i\in[1,k]$.
Hence $w(\cha_{s_1}+\do+\cha_{s_k})$ is an $\NN$-linear combination of elements in $\che R^+$ and is also
equal to $\cha_{s_1}+\do+\cha_{s_k}$, a contradiction. Thus
$w(\cha_{s_i})\in\che R^+$ for any $i\in[1,k]$. This, combined with 
$w(\cha_{s_1})+\do+w(\cha_{s_k})=\cha_{s_1}+\do+\cha_{s_k}$ forces 
the equality $w(\cha_{s_i})=\cha_{s_{h(i)}}$ for all $i\in[1,k]$ where
$h:[1,k]@>>>[1,k]$ is a permutation. Note that $h$ necessarily commutes with the cyclic permutation
of $[1,k]$ induced by $\io$ hence it is a power of this cyclic permutation. Moreover we have $h^2=1$ hence $h=1$
unless $k=2$.

Assume first that $h=1$. We have 

$w(\cha_{s_1})=\cha_{s_1}$, $(s_1w)(\cha_{s_2}=\cha_{s_2}$, $(s_{k-1}\do s_1w)(\cha_{s_k})=\cha_{s_k}$, 
\nl
hence 
$s_1w=ws_1$, $s_2s_1w=s_1ws_2,\do,$  $s_k\do s_1w=s_{k-1}\do s_1ws_k$. By (a) for $G$ we have
$r_{s_1w}-r_w=\pm\cha_{s_1}$, $r_{s_2s_1w}-r_{s_1w}=\pm\cha_{s_2}$,
$r_{s_k\do s_2s_1w}-r_{s_{k-1}\do s_1w}=\pm\cha_{s_k}$. Taking the sum we obtain
$r_{\s w}-r_w=r_{s_k\do s_2s_1w}-r_w=c_1\cha_{s_1}+\do+c_k\cha_{s_k}$ with $c_1,\do,c_k$ in $\{-1,1\}$.
Since $r_{\s w}-r_w$ is fixed by $\io$, so must be $c_1\cha_{s_1}+\do+ c_k\cha_{s_k}$. It follows that
$c_1=\do=c_k$ so that $r_{\s w}-r_w=\pm(\cha_1+\do+\cha_k)$. We see that (a) holds for $\tG$.

Next we assume that $h\ne1$; then $k=2$ and $w(\cha_{s_1})=\cha_{s_2}$, $w(\cha_{s_2})=\cha_{s_1}$.
It follows that $ws_1w=s_2$. We have $s_1s_2w=s_1ws_1\ne w$, $s_1s_2w=s_2ws_2\ne w$. By 1.9 for $G$ we have
$r_{s_1s_2w}=r_{s_1ws_1}=s_1(r_w)$, $r_{s_1s_2w}=r_{s_2ws_2}=s_2(r_w)$.
In particular we have $s_1(r_w)=s_2(r_w)$ hence $(s_1s_2)r_w=r_w$.
We have $s_1(r_w)-r_w\in\ZZ\cha_{s_1}$, $s_2(r_w)-r_w\in\ZZ\cha_{s_2}$. Hence
$r_{s_1s_2w}-r_w\in(\ZZ\cha_{s_1})\cap(\ZZ\cha_{s_2})$. 
We have $(\ZZ\cha_{s_1})\cap(\ZZ\cha_{s_2})=0$ hence  $r_{s_1s_2w}=r'_w$ that is $r_{\s w}=r_w$. We see that
(a) holds for $\tG$. This completes the proof of (a).

\subhead 1.16. Proof of Theorem 0.2\endsubhead
The map $W_2@>>>L$, $w\m r_w$ in 1.8 satisfies 0.2(i) by definition, satisfies 0.2(ii) and 0.2(iv) by 1.9 and
satisfies 0.2(iii) by 1.15. It satisfies 0.2(v) since $r_w\in\YY_w$ and $w$ acts as multiplication by $-1$ on
$\YY_w$. This proves the existence part of 0.2.

Assume now that $w\m r'_w$ is a map $W_2@>>>L$ satisfying conditions like 0.2(i)-(iii).
We show that $r'_w=r_w$ for $w\in W_2$ by induction on $|w|$.
When $|w|\le1$ this follows from 0.2(i). Now assume that $|w|\ge2$.
Assume first that there exists $s\in S$ such that $|sw|<|w|$ and $sw\ne ws$. By the induction hypothesis we
have $r'_{sws}=r_{sws}$ hence, by 0.2(ii), $s(r'_w)=s(r_w)$ so that $r'_w=r_w$.
Assume next that no such $s$ exists. Then by 1.4, $w$ is the longest element in a standard parabolic
subgroup $W_J$ of $W$ whose center contains $w$. Since $|w|\ge2$ we can find two distinct elements
$s_1,s_2$ of $S$ which are contained in $W_J$. Then $s_1w\in W_2,s_2w\in W_2$ and $|s_1w|<|w|,|s_2w|<|w|$, so that
by the induction hypothesis we have $r'_{s_1w}=r_{s_1w}$, $r'_{s_2w}=r_{s_2w}$.
Now let $s\in S$. If $s\ne s_1$, the coefficient of $\cha_s$ in $r'_w$ is equal to 
the coefficient of $\cha_s$ in $r_w$ (they are both equal to the coefficient of $\cha_s$ in 
$r'_{s_1w}=r_{s_1w}$ (see 0.2(iii)). If $s=s_1$, then $s\ne s_2$ and the coefficient of $\cha_s$ in $r'_w$ is 
equal to the coefficient of $\cha_s$ in $r_w$ (they are both equal to the coefficient of $\cha_s$ in 
$r'_{s_2w}=r_{s_2w}$ (see 0.2(iii)). Thus $r'_w=r_w$. This completes the induction. Theorem 0.2 is proved.

\subhead 1.17\endsubhead
For $w\in W_2,s\in S$ such that $sw=ws$ we define a number $(w:s)\in\{-1,0,1\}$ as follows.
Assume first that $G$ is almost simple, simply laced. 
The root system $\che R_w, R_w$ is simply laced and has no component of type $A_l,l>1$. 
Moreover we have $\cha_s\in\che\Pi_w$. 

If the component containing $\cha_s$ is not of type $A_1$, there is a unique sequence
$\cha_1,\cha_2,\do,\cha_m$ in $\che\Pi_w$ such that 
$\cha_i,\cha_{i+1}$ are joined in the Dynkin diagram of $\che R_w$
for $i=1,2,\do,m-1$, $\cha_1=\cha_s$, and $\cha_m$ corresponds to a branch point of the Dynkin diagram of $\che R_w$;
if the component containing $\cha_s$ is of type $A_1$ we define $\cha_1,\cha_2,\do,\cha_m$ as the sequence with one term $\cha_s$ (so that
$m=1$. We define $(w:s)=(-1)^m$ if $|sw|<|w|$ and $(w:s)=(-1)^{m+1}$ if $|sw|>|w|$.
Next we assume that $G$ is almost simple, not simply laced.
Then $G$ can be regarded as a fixed point set of an automorphism  of a simply connected almost simple, simply laced group $G'$
(as in 1.14) with Weyl group $W'$, a Coxeter group with a length preserving automorphism $W'@>>>W'$ with
fixed point set $W$. When $s$ is regarded as an element of $W'$, it is a product of $k$ commuting simple reflections
$s'_1,s'_2,\do,s'_k$ of $W'$; here $k\in\{1,2,3\}$. If $k\in\{0,3\}$ then we define $(w:s)$ for $W$ to be $(w:s_i)$ for $G'$ where
$i$ is any element of $\{1,2,3\}$. If $k=2$ we have either $ws_1=s_1w$, $ws_2=s_2w$
(and $(w:s)$ for $G$ is defined to be $(w:s_1)=(w:s_2)$ for $G'$) or $ws_1=s_2w$, $ws_2=s_1w$ (and $(w:s)$ for $G$ is defined to be $0$.)
We now drop the assumption that $G$ is almost simple. Let $G''$ be the almost simple factor of $G_{der}$ with Weyl group $W''\sub W$
such that $s\in W''$ and let $w''$ be the $W'$-component of $w$. Then $(w:s)$ for $G$ is defined to be $(w'':s)$ for $G''$ (which is
is defined as above).

\mpb

The proof of  1.15(a) yields the following refinement of 1.15(a).

\proclaim{Lemma 1.18} Let $w\in W_2,s\in S$ be such that $sw=ws$. We have 

(a) $r_{sw}=r_w+(w:s)\cha_s$. 
\endproclaim

\head 2. The elements $b_w$\endhead
Assume that we are in the setup of 0.1. We have the following result.

\proclaim{Lemma 2.1} (a) Let $W_J$ be the parabolic subgroup of $W$ generated by $J\sub S$; we assume that 
$W_J$ is an irreducible Weyl group and that the centre of $W_J$ contains the longest element $w_J$ of $W_J$. Let 
$\a=\a_J\in R$ be the unique root such that $\a=\sum_{s\in J}u_s\a_s$ with $u_s\in\NN$ and $\sum_su_s$ as large as 
possible. We have $\ds_\a^2=r_{s_\a}(\e)$.

(b) We have $\dw_J^2=r_{w_J}(\e)$.
\endproclaim
We can assume that $W$ is irreducible. We denote the simple roots by $\{\a_i;i\in[1,l]\}$ as in 
\cite{\BOU} and the corresponding simple reflections as $\{s_i;i\in[1,l]\}$. 
We write $\e_i$ instead of $\e_{s_i}$. We write $i_1i_2\do i_k$ 
instead of $\ds_{i_1}\ds_{i_2}\do\ds_{i_k}$.

We  note that (a) does not hold in general if $w_J$ is not central in $W_J$.
For example if $W=W_J$ is of type $A_2$ we have $(121)^2=121121=12\e_121=1\e_11=\e_1=\text{unit element}$
and $r_{121}(\e)=\e_1\e_2$.

We prove (a). By an argument in the proof of 1.14(a), we can reduce the general case to the case where
$G$ is simply laced. Moreover, we can assume that $J=S$ hence $W_J=W$.
In this case the proof of (a) is case by case. 

Type $A_1$. We have $\ds_\a^2=\ds_1^2=\e_1=r_\a(\e)$.

Type $D_l$, $l=2n\ge4$. We have
$$\ds_\a=234\do(l-2)(l-1)(l)(l-2)\do212\do(l-2)(l)(l-1)(l-2)\do 432.$$
A direct computation shows that $\ds_\a^2=\e_{l-1}^n\e_l^n$ and this is also equal to $r_{s_\a}(\e)$ (see 1.2). 
For example if $l=4$ we have
$$\align&234212342234212342=23421234\e_234212342=2342123\e_2\e_43212342\\&=
234212\e_2\e_3\e_4212342=23421\e_3\e_412342=2342\e_3\e_42342\\&=
234\e_2\e_3\e_4342=23\e_2\e_332=2\e_22=\text{unit element}.\endalign$$
Type $E_7$. We have 
$$\ds_\a=134567243156432545234651342765431.$$
(See \cite{\SQINT}.) A direct computation (as for $D_4$ above) shows that $\ds_\a^2=\e_3\e_5\e_7=r_{\ds_\a}(\e)$.
(See 1.2.)

Type $E_8$. We have
$$\ds_\a=876542314563457624587634524313425436785426754365413245678.$$
(See \cite{\SQINT}.) A direct computation (as for $D_4$ above) shows that $\ds_\a^2=\e_2\e_5\e_7=r_{\ds_\a}(\e)$.
(See 1.2.) This proves (a).

We prove (b). Let $w$ be the longest element of $W$. By 1.8(a) we have $w=\prod_{\b\in\che\ce}s_\b$ with $s_\b$ 
commuting with each other; moreover each $s_\b$ is of the form $s_{\a_{J'}}$ where $J'$ is like $J$ in the lemma
hence (a) is applicable to it. Thus $\ds_\b^2=r_{s_\b}(\e)$. From the description of $\che\ce$ in 1.1
we see that $|w|=\sum_{\b\in\che\ce}|s_\b|$ hence $\dw=\prod_{\b\in\che\ce}\ds_\b$ and 
$\dw^2=\prod_{\b\in\che\ce}\ds_\b^2$ (using the fact the $s_\b$ commute). Using (a) for $s_\b$ we obtain
$\dw^2=\prod_{\b\in\che\ce}r_{\ds_\b}(\e)$ hence $\dw^2=r_w(\e)$. The lemma is proved.

\subhead 2.2\endsubhead
In this subsection we prove the following weak version of Theorem 0.3.

(a) {\it For any $w\in W_2$ one can find $b_w\in L/2L$ such that 
$b_w(\e)w(b_w(\e))=r_w(\e)\dw^2$ or equivalently $(\dw b_w(\e))^2=r_w(\e)$.}
\nl
We argue by induction on $|w|$. If $|w|=0$ we can take $b_w=0$. Now assume that $|w|\ge1$. Assume first that 
there exists $s\in S$ such that $|sw|<|w|$ and $sw\ne ws$. Then $|sws|=|w|-2$. Using the induction hypothesis 
applied to $w':=sws$ and 0.2(ii) we see that we can find $b\in L/2L$
such that $b(\e)w'(b(\e))=(s(r_w))(\e)\dw'{}^2$. Let $b'=s(b)+\cha_s\in L/2L$ so that
$b(\e)=s(b'(\e))\e_s$
We have $s(b'(\e)\e_s)s(w(b'(\e)\e_s))=s(r_w(\e))\dw'{}^2$
hence $b'(\e)\e_sw(b'(\e))w(\e_s)=r_w(\e)s(\dw'{}^2)$. We show that $b'(\e)w(b'(\e))=r_w(\e)\dw^2$. 
It is enough to show that $\e_sw(\e_s)\dw^2=s(\dw'{}^2)$ or that 
$\e_sw(\e_s)\ds\dw'\ds\ds\dw'\ds=\ds\i\dw'{}^2\ds$ or that $\e_sw(\e_s)\ds\dw'\e_s=\ds\i\dw'$
or that $\e_sw(\e_s)\ds w'(\e_s)=\e_s\ds$.
This is immediate. Thus we can take $b_w=b'$ and (a) holds for $w$.

Next we assume that no $s\in S$ as above can be found. Then, by Lemma 1.4, $w$ is the 
longest element in a standard parabolic subgroup of $W$ whose centre contains $w$. By 2.1 we have 
$\dw^2r_w(\e)=1$. Thus we can take $b_w=0$. This completes the proof of (a).

Note that the elements $b_w$ do not necessarily satisfy conditions 0.3(ii),(iii). The interest in
proving the weaker result (a) is that unlike the proof of 0.3, it does not rely on computer calculations.

\subhead 2.3\endsubhead
We prove the uniqueness statement in Theorem 0.3. The argument is similar to that in the proof of uniqueness in 0.2.
Assume that $b',b''$ are two functions $W_2@>>>L/2L$ satisfying conditions like (i),(ii),(iii) in 0.3. We 
show that $b'(w)=b''(w)$ for $w\in W_2$ by induction on $|w|$. 
When $|w|\le1$ this follows from 0.3(i). Now assume that $|w|\ge2$.
Assume first that there exists $s\in S$ such that $|sw|<|w|$ and $sw\ne ws$. By the induction hypothesis we
have $b'_{sws}=b''_{sws}$ hence, by 0.3(ii), $s(b'_w)+\cha_s=s(b''_w)+\cha_s$ so that $b'_w=b''_w$.
Assume next that no such $s$ exists. Then by 1.4, $w$ is the longest element in a standard parabolic
subgroup $W_J$ of $W$ whose center contains $w$. Since $|w|\ge2$ we can find two distinct elements
$s_1,s_2$ of $S$ which are contained in $W_J$. Then $s_1w\in W_2,s_2w\in W_2$ and $|s_1w|<|w|,|s_2w|<|w|$, so that
by the induction hypothesis we have $b'_{s_1w}=b''_{s_1w}$, $b'_{s_2w}=b''_{s_2w}$.
Now let $s\in S$. If $s\ne s_1$, the coefficient of $\cha_s$ in $b'_w$ is equal to 
the coefficient of $\cha_s$ in $b''_w$ (they are both equal to the coefficient of $\cha_s$ in 
$b'_{s_1w}=b''_{s_1w}$ (see 0.3(iii)). If $s=s_1$ then $s\ne s_2$ and the coefficient of $\cha_s$ in $b'_w$ is 
equal to the coefficient of $\cha_s$ in $b''_w$ (they are both equal to the coefficient of $\cha_s$ in 
$b'_{s_2w}=b''_{s_2w}$ (see 0.3(iii)). Thus $b'_w=b''_w$. This completes the inductive proof of uniqueness.

\subhead 2.4\endsubhead
We sketch a proof of the existence part of Theorem 0.3 in the setup of 1.10. In this case the set 
$\Sigma$ of simple coroots consists of
$$e_1-e_2,e_2-e_3,\do,e_{l-1}-e_l,e_{l-1}+e_l.$$
Let $w\in W_2$. For any two element subset $\{\b,\b'\}$ of $\che\ce_w$ we define a subset 
$\cm_{\b,\b'}\sub\Sigma$ as follows.

(a) Assume that $\{\b,\b'\}=\{e_i-e_j,e_k-e_h\}$ where $i<j,k<h$, $i\ne k,i\ne h,j\ne k,j\ne h$.
Then $\cm_{\b,\b'}$ consists of all $e_a-e_{a+1}\in\Sigma$ such that $i\le a<a+1\le j$, $k\le a<a+1\le h$.

(b) Assume that $\{\b,\b'\}=\{e_i-e_j,e_k+e_h\}$ where $i<j,k<h$, $i\ne k,i\ne h,j\ne k,j\ne h$.
 Then $\cm_{\b,\b'}$ consists of all $e_a-e_{a+1}\in\Sigma$ such that $i\le a<a+1\le j$, $k\le a<a+1\le h$.

(c) Assume that $\{\b,\b'\}=\{e_i+e_j,e_k+e_h\}$ where $i<j,k<h$, $i\ne k,i\ne h,j\ne k,j\ne h$.
 Then $\cm_{\b,\b'}$ consists of all $e_a-e_{a+1}\in\Sigma$ such that 
$a=n-2\mod2$, $i\le a<a+1\le j$, $k\le a<a+1\le h$.

(d) Assume that $\{\b,\b'\}=\{e_i-e_j,e_i+e_j\}$ where $i<j$.  Then $\cm_{\b,\b'}$ consists of all 
$e_a-e_{a+1}\in\Sigma$ such that $a=n-2\mod2$, $i\le a<a+1\le j$.
\nl
Let $S'$ be a halving of $S$. We have $r_w=\sum_{s\in S}c_s\cha_s$
where $c_s\in\NN$. We set $r_w^{S'}=\sum_{s\in S}c'_s\cha_s\in L/2L$ where $c'_s=c_s$ if $s\in S'$, $c'_s=0$
if $s\in S-S'$. We define
$$b_w=r_w^{S'}+\sum_{\b,\b'}\sum_{s\in S;\cha_s\in \cm_{\b,b'}}\cha_s$$
where $\{\b,\b'\}$ runs through all 2 element subsets of $\che\ce_w$.
One can verify that the elements $b_w,w\in W_2$ satisfy conditions (i)-(iv) in 0.3.
This, together with 2.3 proves Theorem 0.3 whe $G$ is almost simple of type $D_l$ (except for condition 0.3(v)).

Now if $G$ is adjoint of type $A_l$, $l\ge1$, then $G$ can be regarded as the adjoint group of a Levi subgroup 
of a parabolic subgroup in a group of type $D_{l'}$ for some $l'$ such that $l<l'\ge4$ and Theorem 0.3 for $G$
can be deduced from the results above (for type $D_{l'}$) (except for condition 0.3(v)).

Next, the argument in the proof of uniqueness in 2.3 can be viewed as an inductive method to compute $b_w$ in 0.3
for any $w\in W_2$ by induction on $|w|$. This can be used to prove the existence statement in 0.3(i)-(iv)
in any given case with a powerful enough computer. We have used this method to prove the 
existence statement in 0.3(i)-(iv) for $G$ of type $E_8$. (I thank Gongqin Li for carrying out the
programming in GAP using the CHEVIE package.) Then 
0.3(i)-(iv) automatically holds for $G$ of type $E_7$ and $E_6$.
We see that Theorem 0.3 holds for any simply laced $G$ (except for condition 0.3(v)).

\subhead 2.5\endsubhead
We show:

(a) {\it If 0.3(i)-(iv) is assumed to hold for $G$ then 0.3(v) holds for $G$.}
\nl
We prove the equality in 0.3(v) for $w\in W_2$ by induction on $|w|$. If $|w|=0$ the result is obvious.
Assume first that there exists $s\in S$ such that $sw\ne ws$, $|sw|<|w|$. We have $|sws|=|w|-2$.
By the induction hypothesis we have $w'(b_{w'}(\e))b_{w'}(\e)=r_{w'}(\e)\dw'{}^2$ where $w'=sws$.
As in the proof in 2.2 we deduce that $b':=s(b_{w'})+\cha_s$ satisfies $w(b'(\e))b'(\e)=r_w(\e)\dw^2$.
By 0.3(ii) we have $b'=b_w$. Thus 0.3(v) holds for $w$.
Next we assume that no $s$ as above can be found. Then, by 1.4, $w$ is the longest element of a standard 
parabolic subgroup $W_J$ of $W$ and $w$ is in the centre of $W_J$. In this case, using 2.1, we see that it is 
enough to show that $w(b_w(\e))=b_w(\e)$. From the definition we see that $b_w=\sum_{s\in J}a_s\cha_s$ where 
$a_s\in\{0,1\}$. Hence to show that
$w(b_w(\e))=b_w(\e)$ it is enough to show that for any $s\in J$ we have $w(\cha_s)=\cha_s$ in $L/2L$. This is clear
since $w(\cha_s)=-\cha_s$ in $L$. This completes the proof of (a).

We see that Theorem 0.3 holds for any simply laced $G$.

\subhead 2.6\endsubhead
In this subsection we assume that $G,\io,\tG,\tT,\tW,\tS,\tY,||_\io,\tW_2$ 
are as in the proof of 1.14(a). Assume that $S'$ is a halving of $S$ such that $\io(S')=S'$. (Such a halving 
exists.) Assume also that $w\m b_w$ is a function $W_2@>>>L/2L$ satisfying 0.3(i)-(iv) for $G$. Let $\tS'$ 
be the subset of $\tS$ consisting of the elements $\s=\prod_ss$ where $s$ runs over an $\io$-orbit in $S'$. 
Clearly, $\tS'$ is a halving of $\tS$ (and any halving of $\tS$ is of this form).
Let $\tL$ be the subgroup of $L$ generated by the coroots of $\tG$. We have canonically
$\tL/2\tL=(L/2L)^\io$. We define a function $\tb:\tW_2@>>>\tL/2\tL$ by $w\m\tb_w=b_w$. (Note that if
$w\in \tW_2$ then $\io(b_w)=b_w$, by the uniqueness statement in 0.3.) We show:

(a) {\it The function $\tW_2@>>>\tL/2\tL$, $w\m\tb_w$ satisfies 0.3(i)-(iv) for $\tG$.}
\nl
We have $\tb_1=0$. Let $\s\in\tS$. Now $\s$ corresponds to an $\io$-orbit $\co$ in $S$. We can view
$\s$ as an element of $W_2$ and we have $\tb_\s=b_\s=\sum_{s\in\co}\cha_s$ if $\co\sub S'$
and $\tb_\s=0$ if $\co\sub S-S'$ (this follows from 0.3(i) for $b$ applied to
each $s\in\co$ and from 0.3(iii) for $G$ applied to $\s\in W_2$). Thus 
0.3(i) holds for $\tb$.

Now let $w\in\tW_2$ and $\s\in\tS$ be such that $\s w\ne w\s$. We write $\s=\{s_1,\do,s_k\}\sub S$.
Then for some $i\in[1,k]$ we have $s_iw\ne ws_i$. Hence for all $i\in[1,k]$ we have $s_iw\ne ws_i$. 
Hence we have $s_1w\ne ws_1, s_2s_1w\ne ws_1s_2,\do,s_k\do s_1w\ne ws_1\do s_k$. By 0.3(ii) for $b$ 
we have $b_w=s_1(b_{s_1ws_1})+\cha_{s_1}$, $s_1(b_{s_1ws_1})=s_1s_2(b_{s_2s_1ws_1s_2})+\cha_{s_2}$,
$s_1\do s_k(b_{s_k\do s_1ws_1\do s_k})=s_1\do s_{k-1}(b_{s_{k-1}\do s_1ws_1\do s_{k-1}})+\cha_{s_k}$
so that $b_{\s w\s}=\s(b_w)+\cha_{s_1}+\do+\cha{s_k}$
that is $\tb_{\s w\s}+\s(\tb_w)+\cha_{s_1}+\do+\cha{s_k}$. We see that 0.3(ii) holds for $\tb$.

Next, let $w\in\tW_2\cap\tW$, $\s\in\tS$ be such that $\s w=w\s$.
We write $\co=\{s_1,\do,s_k\}\sub S$, $\s=s_1\do s_k$
where $\io$ permutes $s_1,s_2,\do,s_k$ cyclically: $s_1\m s_2\m\do s_k\m s_1$. (Note that $k\le 3$.)
As in the proof of 1.15 we have $w(\cha_{s_i})=\cha_{s_{h(i)}}$ for all $i\in[1,k]$ 
where $h$ is a permutation of $[1,k]$ such that $h=1$ unless $k=2$.

Assume first that $h=1$. We have 

$w(\cha_{s_1})=\cha_{s_1}$, $(s_1w)(\cha_{s_2}=\cha_{s_2}$, $(s_{k-1}\do s_1w)(\cha_{s_k})=\cha_{s_k}$
\nl
 hence 
$s_1w=ws_1$, $s_2s_1w=s_1ws_2,\do,$  $s_k\do s_1w=s_{k-1}\do s_1ws_k$,
$|s_1w|>|w|$, $|s_2s_1w|>|s_1w|,\do,|s_k\do s_1w|>|s_{k-1}\do s_1w|$. By 0.3(iii) for $b$ we have
$b_{s_1w}-b_w=l_1\cha_{s_1}$, $b_{s_2s_1w}-b_{s_1w}=l_2\cha_{s_2}$,
$b_{s_k\do s_2s_1w}-b_{s_{k-1}\do s_1w}=l_k\cha_{s_k}$ with $l_1,\do,l_l$ in $\{0,1\}$. Taking the sum we obtain
$b_{\s w}-b_w=b_{s_k\do s_2s_1w}-b_w=l_1\cha_{s_1}+\do+l_k\cha_{s_k}$.
Since $b_{\s w}-b_w$ is fixed by $\io$, so must be $l_1\cha_{s_1}+\do+l_k\cha_{s_k}$. It follows that
$l_1=\do=l_k$ so that $b_{\s w}-b_w=l_1(\cha_1+\do+\cha_k)$. We see that 0.3(iii) holds for $\tb$.
Using repeatedly 0.3(iv) for $b$ we have $\s(b_w)=s_1s_2\do s_k(b_w)=b_w$ (here we use that $G$ is simply 
laced). We see that 0.3(iv) holds for $\tb$ (in this case we have $r_{\s w}-r_w=\pm(\cha_1+\do+\cha_k)$ by 
the proof of 1.15).

Next we assume that $h\ne1$; then $k=2$ and $w(\cha_{s_1})=\cha_{s_2}$, $w(\cha_{s_2})=\cha_{s_1}$.
It follows that $ws_1w=s_2$. We have $s_1s_2w=s_1ws_1\ne w$, $s_1s_2w=s_2ws_2\ne w$. By 0.3(ii) for $b$ 
we have

$b_{s_1s_2w}=b_{s_1ws_1}=s_1(b_w)+\cha_{s_1}=b_w+c_1\cha_1$, 

$b_{s_1s_2w}=b_{s_2ws_2}=s_2(b_w)+\cha_{s_2}=b_w+\c_2\cha_2$ where $c_1,c_2\in\{0,1\}$. 
\nl
It follows that $c_1\cha_1=c_2\cha_2$ in $L/2L$ hence $c_1=c_2=0$
and $b_{s_1s_2w}=b_w$ that is $b_{\s w}=b_w$. We see that 0.3(iii) holds for $\tb$.
By 0.3(ii) for $b$ we have $b_{ws_1s_2}=b_{s_1ws_1}=s_1(b_w)+\cha_1$ and
$b_{ws_1s_2}=b_{s_2ws_2}=s_2(b_w)+\cha_2$ hence $s_1(b_w)+\cha_1=s_2(b_w)+\cha_2$. Applying $s_1$ we obtain
$s_1s_2(b_w)+\cha_2=b_w+\cha_1$ hence $\s(b_w)=b_w+\cha_1+\cha_2$. We see that 0.3(iv) holds for $\tb$
(in this case we have $r_{\s w}-r_w=0$ by the proof of 1.15). This completes the proof of (a).

\subhead 2.7\endsubhead
From 2.6 we see that Theorem 0.3(i)-(iv) can be reduced to the case where $G$ is simply laced.
Using this and the results in 2.5 we see that Theorem 0.3 holds in the general case.

\subhead 2.8\endsubhead
Let $S'$ be a halving of $S$. Then clearly $S-S'$ is a halving of $S$. We define $W_2@>>>L/2L$ by $w\m b^*_w=b_w+r_w$
where $b_w=b_w^{S'}$. We have

(a) {\it $b_w^{S'-S}=b_w^*$.}
\nl
The fact that $b^*$ satisfies 0.3(i) for $S-S'$ is immediate; that it satisfies 0.3(ii) for $S-S'$ 
follows from 1.9(a); that it satisfies 0.3(iii) for $S-S'$ follows from 1.13(a).
Then (a) follows from the uniqueness in 0.3.

\subhead 2.9. Proof of Theorem 0.5\endsubhead
In this subsection we assume that we are in the setup of 0.4. Let $w,c,S'$ be as in 0.5. We must show that
$$\ph(n_{w,c,S'})n_{w,-\ph(c),S'})=1.$$
We write $b_w$ instead of $b_w^{S'}$. We have $\ph(\dw)=\dw$, $\ph(r_w(c))=r_w(\ph(c))$, $\ph(b_w(\e))=b_w(\e)$, 
hence
$$\align&\ph(n_{w,c,S'})n_{w,-\ph(c),S'}=\ph(\dw r_w(c)b_w(\e))\dw r_w(-\ph(c))b_w(\e)\\&
=\dw r_w(\ph(c))b_w(\e)\dw r_w(-\ph(c))b_w(\e)=
\dw^2 w(r_w(\ph(c)))w(b_w(\e))r_w(-\ph(c))b_w(\e)\\&=
w(r_w(\ph(c)))r_w(-\ph(c))r_w(\e)=w(r_w(\ph(c)))r_w(\ph(c)).\endalign$$
This equals $1$ since $w(r_w(\ph(c))=r_w(\ph(c)\i)$ by 0.2(v). Theorem 0.5 is proved.

\subhead 2.10\endsubhead
Assume now that $G$ is almost simple. Then there are exactly
two halvings $S$, $S'-S$ for $S$. Let $w\in W_2$. We note that the family of elements
$\{n_{w,c,S'};c\in\kk^*\}$ coincides with the family of elements 
$\{n_{w,c,S-S'};c\in\kk^*\}$.
Indeed, by 2.8(a) we have 
$$n_{w,c,S-S'}=\dw r_w(c)b^{S'}_w(\e)r_w(\e)=n_{w,\e c,S}.$$

\subhead 2.11\endsubhead
In this subsection we assume that $\kk,G,\ph,\ph',F_q$ are as in 0.4 (in case 0.4(i)). Now

(a) $g_1:g\m g_1g\ph(g_1)\i$ 
\nl
defines an action of $G^{\ph^2}=G(F_{q^2})$ on $G^{\ph'}$.
Indeed for $g_1\in G^{\ph^2}$, $g\in G^{\ph'}$. We have 
$$\align&\ph(g_1g\ph(g_1)\i)g_1g\ph(g_1)\i=\ph(g_1)\ph(g)g_1\i g_1g\ph(g_1)\i\\&=\ph(g_1)\ph(g)g\ph(g_1)\i=
\ph(g_1)\ph(g_1)\i=1\endalign$$ 
and our claim follows. We have $1\in G^{\ph'}$ and the stabilizer of $1$ for the
action above is $G^\ph$. Thus we have an injective map $G^{\ph^2}/G^\ph@>>>G^{\ph'}$.
We show that this is a bijection. Let $g\in G^{\ph'}$. By Lang's theorem we have $g=g_1\ph(g_1)\i$
for some $g_1\in G$. We have $g\ph(g)=1$ hence
$g_1\ph(g_1)\i\ph(g_1)\ph^2(g_1)\i=1$ that is $g_1\ph^2(g_1)\i=1$ so that $g_1\in G^{\ph^2}$. We see that
$g$ is in the $G^{\ph^2}$-orbit of $1$. Thus we have the following result.

(b) {\it The action (a) of $G^{\ph^2}=G(F_{q^2})$ on $G^{\ph'}$ is transitive; the stabilizer of $1$ for 
this action is $G^\ph$. Hence $\sha(G^{\ph'})=\sha(G^{\ph^2})/\sha(G^\ph)$.}

\enddocument